\documentclass[12pt,reqno]{amsart}

\usepackage{amscd,amssymb}

\newcommand{\Z}{{\mathbb Z}}

\newcommand{\C}{{\mathbb C}}
\newcommand{\R}{{\mathbb R}}
\renewcommand{\P}{{\mathbb P}}

\newcommand{\AAA}{{\mathcal A}}
\newcommand{\BB}{{\mathcal B}}
\newcommand{\CC}{{\mathcal C}}

\newcommand{\OO}{{\mathcal O}}
\newcommand{\PP}{{\mathcal P}}

\newcommand{\TT}{{\mathcal T}}

\newcommand{\www}{\widetilde}

\newcommand{\paa}{\partial}

\newcommand{\nnn}{\nabla}

\DeclareMathOperator{\rank}{rank}

\DeclareMathOperator{\supp}{supp}

\theoremstyle{plain}
\newtheorem{lemma}{Lemma}[section]
\newtheorem{theorem}[lemma]{Theorem}

\theoremstyle{definition}
\newtheorem{definition}[lemma]{Definition}

\newtheorem{remark}[lemma]{Remark}
\newtheorem{remarks}[lemma]{Remarks}

\newtheorem{notations}[lemma]{Notations}

\begin{document}

\title[Potentials of a Frobenius like structure]
{Potentials of a Frobenius like structure\\ and $m$ bases of a vector space} 

\author[C. Hertling and A. Varchenko]{Claus Hertling and Alexander Varchenko}

\address{Lehrstuhl f\"ur Mathematik VI,
Universit\"at Mannheim, A5,6, 68131 Mannheim, Germany}
\email{hertling\char64 math.uni-mannheim.de}

\address{Department of Mathematics, University of North Carolina 
at Chapel Hill, Chapel Hill, NC 27599-3250, USA}

\email{anv@email.unc.edu}

%\date{August 18, 2016}



\thanks{This work was supported by the DFG grant He2287/4-1
(SISYPH), the second author was supported in part by NSF grant DMS-1362924
          and Simons Foundation grant \#336826}

\maketitle

%\clearpage

\begin{abstract}
This paper proves the existence of potentials of the first
        and second kind of a Frobenius like structure in a frame
        which encompasses families of arrangements.

       Surprisingly the proof is based on the study of finite sets of vectors
     in a finite-dimensional vector space $V$.  Given a natural number $m$ 
      and a finite set  $(v_i)$ of vectors we give a necessary and sufficient
      condition to find in the set $(v_i)$ $m$ bases of  $V$.  If $m$ bases in $(v_i)$
      can be selected,  we define elementary transformations of such
      a selection and show that any two selections are
      connected by a sequence of elementary transformations. 
\end{abstract}

%\clearpage

%\tableofcontents

%\clearpage

% section 1
\section{Introduction and main results}\label{s1}
\setcounter{equation}{0}

\noindent
A Frobenius manifold comes equipped locally with a potential.
If one gives a definition which does not mention this potential
explicitly, one nevertheless obtains it immediately by the following 
elementary fact: Let $z_i$ be the coordinates on $\C^n$
and $\paa_i=\frac{\paa}{\paa z_i}$ be the coordinate vector fields.
Let $M$ be a convex open subset of $\C^n$
and $\TT_M$ be the holomorphic tangent bundle of $M$.
Let $A:\TT_M^3\to\OO_M$ be a symmetric map such that also
$\paa_i A(\paa_j,\paa_k,\paa_l)$ is symmetric in $i,j,k,l$.
Then a potential $F\in\OO_M$ with 
$\paa_i\paa_j\paa_k F=A(\paa_i,\paa_j,\paa_k)$ exists. On Frobenius manifolds see \cite{D1, D2, M}.

This paper is devoted to a nontrivial generalization of this fact.
The generalization turns up in the theory of families of arrangements as in \cite[ch. 3]{V2}. 
Theorem \ref{t1.2} below gives the main result.
Definition \ref{t1.1} prepares the frame and the used notions.

\begin{definition}\label{t1.1}
(a) A {\it Frobenius like structure of order} $(n,k,m)\in\Z_{>0}^3$ 
with $n\geq k$ is a tuple $(M,K,\nnn^K,C,S,\zeta,V,(v_1,...,v_n))$ 
with the following properties.
$M$ is an open subset of $\C^n$, $K\to M$ is a holomorphic 
vector bundle on $M$ with flat connection $\nnn^K$, 
$C$ is a Higgs field on $K$ with $\nnn^K(C)=0$, 
$S$ is a $\nnn^K$-flat $m$-linear form $S:\OO(K)^m\to\OO_M$,
which is Higgs field invariant, i.e. 
\begin{eqnarray}\label{1.1}
S(C_Xs_1,s_2,...,s_m)=S(s_1,C_Xs_2,...,s_m)=...=S(s_1,s_2,...,C_Xs_m)
\end{eqnarray}
for $s_1,s_2,...,s_m\in\OO(K)$ and $X\in \TT_M$,
$\zeta$ is a global section in $K$,
$V$ is a $k$-dimensional $\C$-vector space,
and $(v_1,...,v_n)$ is a tuple of vectors in $V$ 
with $\langle v_1,...,v_n\rangle=V$
such that the following holds. Denote $J:=\{1,...,n\}$.
For any $\{i_1,...,i_k\}\subset J$ with $(v_{i_1},...,v_{i_k})$ 
a basis of $V$, the section $C_{\paa_{i_1}}...C_{\paa_{i_k}}\zeta$
is $\nnn^K$-flat.

\medskip
(b) Some notations: In the situation of (a), a subset $I=\{i_1,...,i_k\}\subset J$
is called {\it independent} if $(v_{i_1},...,v_{i_k})$ is a basis of $V$.
Then the differential operator $\paa_I:=\paa_{i_1}...\paa_{i_k}$ and 
the endomorphism $C_I:=C_{\paa_{i_1}}...C_{\paa_{i_k}}:\OO(K)\to\OO(K)$ 
are well defined (they do not depend on the chosen order 
of the elements $i_1,...,i_k$).
Independent subsets exist because of $\langle v_1,...,v_n\rangle=V$.

\medskip
(c) In the situation of (a), a {\it potential of the first kind}
is a function $Q\in \OO_M$ with
\begin{eqnarray}\label{1.2}
\paa_{I_1}...\paa_{I_m}Q=S(C_{I_1}\zeta,...,C_{I_m}\zeta)
\end{eqnarray}
for any $m$ independent subsets $I_1,...,I_m\subset J$.
A {\it potential of the second kind} is a function $L\in\OO_M$ with 
\begin{eqnarray}\label{1.3}
\paa_i\paa_{I_1}...\paa_{I_m}L=S(C_{\paa_i}C_{I_1}\zeta,...,C_{I_m}\zeta)
\end{eqnarray}
for any $m$ independent subsets $I_1,...,I_m\subset J$ and any $i\in J$.
\end{definition}

Notice that the right-hand side of (\ref{1.2}) does not depend on a point in $M$
    and the right-hand side of (\ref{1.3}) can be interpreted as the matrix
    element of the operator $C_{\paa_i}$.

\begin{theorem}\label{t1.2}
Let $(M,K,\nnn^K,C,S,\zeta,V,(v_1,...,v_n))$ be a Frobenius like 
structure of some order $(n,k,m)\in\Z_{>0}^3$. 
Then locally (i.e. near any $z\in M\subset \C^n$)
potentials of the first and second kind exist.
\end{theorem}

At the end of the paper, several remarks discuss the 
case of arrangements and the relation to Frobenius manifolds
and F-manifolds. But the detailed study of the
case of arrangements is left for the future.
Remark \ref{t4.4} (iii) states some other wishes
for the future. 

The proof of theorem \ref{t1.2} uses a fact from
linear algebra which has (to our knowledge) not been
noticed before.

\begin{theorem}\label{t1.3}
Let $k$ and $m\in\Z_{>0}$. 
Let $V$ be a $k$-dimensional vector space over some field $K$.
Let $(v_1,...,v_{mk})$ be a tuple of vectors in $V$.
It can be split into $m$ bases of $V$ if and only if the following
condition holds: For any vector subspace $U\subset V$
\begin{eqnarray}\label{1.4}
|\{i\in\{1,...,mk\}\, |\, v_i\in U\}|\leq m\cdot\dim U.
\end{eqnarray}
\end{theorem}

The theorem is proved in section \ref{s2}. 
The proof is surprisingly nontrivial. Section \ref{s3}
applies an implication of it to a combinatorial
situation which in turn is needed in the proof of the
main theorem \ref{t1.2} in section \ref{s4}.
Section \ref{s4} concludes with some remarks.

\smallskip
The second author thanks  MPI in Bonn for hospitality during his visit in 2015-2016.

\section{Linear algebra: tuples of vectors giving several bases}\label{s2}
\setcounter{equation}{0}

In this section theorem \ref{t1.3} and some consequences of if 
will be proved. Part (a) of theorem \ref{t2.4} below is a slight 
generalization of theorem \ref{t1.3}. 
Part (b) is a consequence of (a). Only part (b)
will be used later, and only in the case $l=1$.

\noindent
\begin{notations}\label{t2.1}
Let $\AAA\neq\emptyset$ be some set.
An unordererd tuple of elements of $\AAA$
is an element of $\Z_{\geq 0}[\AAA]$.
It is denoted $T=\sum_{a\in\AAA}T(a)\cdot[a]$,
with $T(a)\in\Z_{\geq 0}$ and $|T|:=\sum_{a\in\AAA}T(a)<\infty$.
For simplicity, it is called {\it system} of elements of $\AAA$
instead of {\it unordererd tuple}. If $|T|=k\in\Z_{\geq 0}$,
then it is also called a $k$-system.
Its support is the finite set 
$\supp T:=\{a\in\AAA\, |\, T(a)>0\}\subset\AAA$.
Define for any subset $\BB\subset\AAA$ the number
\begin{eqnarray}\label{2.1}
\mu(T,\BB):=\sum_{a\in \BB\cap\supp T}T(a)\in\Z_{\geq 0}.
\end{eqnarray}
The set $\Z_{\geq 0}[\AAA]$ is a monoid and is contained in the
group $\Z[\AAA]$. The map 
\begin{eqnarray}\label{2.2}
d_H:\Z[\AAA]\times\Z[\AAA]\to\Z_{\geq 0},\quad (T_1,T_2)\mapsto
\sum_{a\in\AAA}|T_1(a)-T_2(a)|
\end{eqnarray}
is a metric on $\Z[\AAA]$. 
On $\Z[\AAA]$ and $\Z_{\geq 0}[\AAA]$ one has the partial ordering
$\leq$ with
\begin{eqnarray}\label{2.3}
S\leq T\iff S(a)\leq T(a)\quad\forall\ a\in\AAA.
\end{eqnarray}
If $S$ and $T$ are systems with $S\leq T$ then $S$ is a 
{\it subsystem} of $T$.
%Any two systems have the maximum $\max(T,S)$ with 
%$\max(T,S)(a):=\max(T(a),T(b))$ and the analogously defined
%minimum $\min(T,S)$. They satisfy $T+S=\max(T,S)+\min(T,S)$.

In the proof of theorem \ref{t2.4} the following notation will
be useful. There an $m\in\Z_{>0}$ will be fixed. Then
for $j\in \Z_{>0}$ let $[j]\in\{1,...,m\}$ be the unique
number with $j\equiv [j]\mod m$.
\end{notations}

\begin{definition}\label{t2.2}
Let $k,m\in\Z_{>0}$, let $l\in\Z_{\geq 0}$, let $K$ be a field, 
let $V$ be a $k$-dimensional $K$-vector space. In the following, 
{\it system} means a system of elements of $V$. 

\begin{list}{}{}
\item[(a)] A system $T$ generates the subspace
$\langle T\rangle :=\langle v\, |\, v\in\supp T\rangle\subset V$.
\item[(b)] 
A system $T$ is a basis of $V$ if $\langle T\rangle=V$ and 
if $|T|=k$ (so the support $\supp T$ is a basis of $V$, 
and all $T(a)\in\{0;1\}$).
\item[(c)]
A {\it strong decomposition} of an $(mk+l)$-system $T$ is a 
decomposition $T=T^{(1)}+...+T^{(m+1)}$ into $m$ $k$-systems 
$T^{(1)},...,T^{(m)}$ and one $l$-system $T^{(m+1)}$ 
such that $T^{(1)},...,T^{(m)}$ are
bases of $V$ (and $T^{(m+1)}$ is an arbitrary $l$-system; 
e.g. if $l=0$ then $T^{(m+1)}=0$ automatically).
\item[(d)] 
An $(mk+l)$-system is {\it strong} if it admits a strong decomposition.
\item[(e)] 
An $(mk+l)$-system $T$ is {\it qualified} if it satisfies the conditions:
\begin{eqnarray}\label{2.4} 
\mu(T,U)\leq l+m\cdot\dim U\quad\forall\ \textup{subspaces }U\subset V.
\end{eqnarray}
\item[(f)]
In the case of an $(mk+l)$-system $T$ with $l\geq 1$, define the subset
\begin{eqnarray}\label{2.5}
A_1(T)&:=&\{a\in\supp T\, |\, \exists\ \textup{a strong decomposition }\\
&&T=T^{(1)}+...+T^{(m+1)}\textup{ with }a\in \supp T^{(m+1)}\} .\nonumber
\end{eqnarray}
Of course, if  $l\geq 1$, then $A_1(T)\neq\emptyset\iff T$ is strong. 
\end{list}
\end{definition}

\begin{lemma}\label{t2.3}
Let $k,m\in\Z_{>0}$, let $l\in\Z_{\geq 0}$, let $K$ be a field, 
let $V$ be a $k$-dimensional $K$-vector space. In the following, 
{\it system} means a system of elements of $V$. 

\begin{list}{}{}
\item[(a)] 
A strong $(mk+l)$-system is qualified, i.e. the conditions \eqref{2.4}
are necessary for being strong.
\item[(b)]
Let $T$ be a qualified $(mk+l)$-system. If $U_1,U_2\subset V$ 
are two subspaces with $\mu(T,U_i)=l+m\dim U_i$, 
then also $\mu(T,U_1\cap U_2)=l+m\dim U_1\cap U_2.$
Therefore there is a unique minimal (with respect to inclusion) 
subspace with this property. The intersection of it with $\supp T$
is called $A_2(T)$. The subspace itself is $\langle A_2(T)\rangle$.
\item[(c)] Let $T$ be a qualified $(mk+l)$-system with $l\geq 1$.
Then $A_1(T)\subset A_2(T).$
\end{list}
\end{lemma}

{\bf Proof:}
(a) Let $T^{(1)}+...+T^{(m+1)}$ be a strong decomposition of $T$.
Each of the bases $T^{(1)},...,T^{(m)}$ contains at most
$\dim U$ elements of $U$. Therefore $\mu(T,U)\leq l+m\dim U$.

(b) The assumption $\mu(T,U_1\cap U_2)\leq (l-1)+m\dim U_1\cap U_2$
leads to a contradiction as follows. If it holds, then 
\begin{eqnarray*}
&&\mu(T,U_1+U_2)\\&=&\mu(T,(U_1+U_2)-(U_1\cup U_2))
+\mu(T,U_1)+\mu(T,U_2)-\mu(T,U_1\cap U_2) \\
&\geq& 0 + (l+m\dim U_1)+(l+m\dim U_2)
-(l-1+m\dim U_1\cap U_2)\\
&=& l+1+m\dim (U_1+U_2),
\end{eqnarray*}
which  contradicts the condition \eqref{2.4} for the qualified 
$(mk+l)$-system $T$. Because the subspace is minimal with the
given property, it is generated by its intersection $A_2(T)$ with $\supp T$.

(c) This follows from (a): If no strong decomposition exists then
$A_1(T)=\emptyset$. If $T^{(1)}+...+T^{(m+1)}$ is a strong 
decomposition of $T$ then each of the bases $T^{(1)},...,T^{(m)}$ 
contains at most $\dim \langle A_2(T)\rangle$ elements of 
$\langle A_2(T)\rangle$. Because of 
$\mu(T,A_2(T))=l+m\dim \langle A_2(T)\rangle$, 
the $l$-system $T^{(m+1)}$ is completely filled with elements of
 $A_2(T)$. \hfill$\Box$ 

\bigskip

\begin{theorem}\label{t2.4}
Let $k,m\in\Z_{>0}$, let $l\in\Z_{\geq 0}$, let $K$ be a field, 
let $V$ be a $k$-dimensional $K$-vector space. In the following, 
{\it system} means a system of elements of $V$. 

\begin{list}{}{}
\item[(a)] 
An $(mk+l)$-system $T$ is strong if and only if it is qualified,
i.e. the conditions \eqref{2.4} are necessary and sufficient for being strong.
\item[(b)]
Let $T$ be a strong $(mk+l)$-system with $l\geq 1$. 
Then $A_1(T)= A_2(T).$
\end{list}
\end{theorem}

Part (a) of the theorem says that the conditions \eqref{2.4} 
are also sufficient for the system to be strong.
As the proof of this fact is surprisingly hard, it is called a theorem.
Part (a) for $l=0$ is theorem \ref{t1.3}.
First the case $l=0$ is proved. The generalization to $l\geq 1$
is an easy consequence. Part (b) of theorem \ref{t2.4}
is an easy consequence of (the proof of) part (a).
It improves part (c) of lemma \ref{t2.3}.

{\bf Proof:}
(a) Because of lemma \ref{t2.3} (a), it rests to show that a qualified
$(mk+l)$-system is strong. Let $T$ be a qualified $(mk+l)$-system.
We can suppose $0\notin\supp T$, i.e. $T(0)=0$. 
If $T(0)>0$, then $T(0)=\mu(T,\{0\})\leq l$ 
by \eqref{2.4}. Then the system $T-T(0)\cdot[0]$ is a qualified
$(mk+(l-(T(0)))$-system.
It is sufficient to prove that this system is strong. Then also $T$ is strong.
Therefore suppose $0\notin\supp T$. Furthermore, 
we can suppose that $T(b)=1$ for all $b\in\supp T$. Because if this does not
hold, one can rescale the vectors in $T=\sum_{j=1}^{mk+l}[b_j]$
with suitable scalars $\lambda_j\in K^*$ so that the new system
$\www T:=\sum_{j=1}^{mk+l}[\lambda_jb_j]$ satisfies $\www T(b)=1$
for all $b\in\supp \www T$. The new system $\www T$ is still qualified.
It is sufficient to prove that $\www T$ is strong.
Then also the old system $T$ is strong.

In a first step, the case $l=0$ will be proved, i.e. theorem \ref{t1.3}.  
In a second step, the cases $l\geq 1$ will be proved inductively on $l$. 
The first step is more difficult than the second step.

\medskip
{\bf First step:}
Let $T$ be a qualified $mk$-system,
so \eqref{2.4} holds with $l=0$. 
As above, we can assume $0\notin\supp T$ and $T(b)=1$ for all $b\in \supp T$.
Let $T=T^{(1)}+...+T^{(m)}$ be a decomposition
of $T$ into $m$ $k$-systems $T^{(1)},...,T^{(m)}$. 
Suppose that it is not a strong decomposition.
We will construct in the following a new decomposition 
$T=S^{(1)}+...+S^{(m)}$ into $m$ $k$-systems $S^{(1)},...,S^{(m)}$
such that 
\begin{eqnarray}\label{2.6}
\sum_{j=1}^m \dim \langle S^{(j)}\rangle
> \sum_{j=1}^m \dim \langle T^{(j)}\rangle.
\end{eqnarray}
Iterating this construction, one arrives at a strong decomposition of $T$.

We can suppose that $T^{(1)}$ is not a basis, so $\dim\langle T^{(1)}\rangle<k$.
For any $j\in\{1,...,m\}$, choose a subset $B^{(j)}\subset\supp T^{(j)}$
such that $(v_j\, |\, j\in B^{(j)})$ is a basis of 
$\langle T^{(j)}\rangle$. 
For any $b\in \langle T^{(j)}\rangle=\langle B^{(j)}\rangle$ 
denote by $\lambda_j(b,c)\in\C$ for $j\in\{1,...,m\}$ 
and $c\in B^{(j)}$ the unique coefficient with
\begin{eqnarray}\label{2.7}
b=\sum_{c\in B^{(j)}}\lambda_j(b,c)\cdot c.
\end{eqnarray}
Define 
\begin{eqnarray}\label{2.8}
R^{(0)}&:=&\supp T^{(1)}-B^{(1)}\neq\emptyset.
\end{eqnarray}
Now, define a sequence $(R^{(j)})_{j=1,...,N}$ of subsets of $\supp T$ with maximal 
$N\in\Z_{>0}\cup\{\infty\}$ in the following way: If $R^{(0)},...,R^{(j-1)}$
for some $j\in\Z_{>0}$ are defined and 
$\langle R^{(j-1)}\rangle\not\subset \langle T^{([j])}\rangle$
(see the notations \ref{t2.1} for $[j]$) then
stop and set $N:=j-1$. But if 
$\langle R^{(j-1)}\rangle\subset \langle T^{([j])}\rangle$
then define
\begin{eqnarray}\label{2.9}
R^{(j)}&:=&\{ c\in B^{([j])}\, |\, \exists\ b\in R^{(j-1)}\textup{ with }
\lambda_{[j]}(b,c)\neq 0\}\\
&=& \textup{the minimal subset }R^{(j)}\subset B^{([j])}\textup{ with }
\langle R^{(j-1)}\rangle \subset\langle R^{(j)}\rangle. \hspace*{1cm}
\label{2.10}
\end{eqnarray}
This defines a unique sequence, a priori with finite or infinite length $N$.
The length satisfies $N\geq 1$, so $R^{(1)}$ exists, 
because $R^{(0)}\subset \langle T^{(1)}\rangle$.

We claim that the length is finite, so $N\in\Z_{>0}$, and prove this indirectly.
So suppose that $N=\infty$.
Because of \eqref{2.10}, there is a $\theta\in\Z_{>0}$ such that 
$\dim\langle R^{(j)}\rangle$ is constant for $j\geq \theta$
and the spaces $\langle R^{(j)}\rangle\subset V$ coincide for all $j\geq\theta$.
Then 
\begin{eqnarray}\label{2.11}
\mu(T,\langle R^{(\theta)}\rangle)
&\geq & |R^{(0)}| + \sum_{j=0}^{m-1}|R^{(\theta+j)}|\\
&=& |R^{(0)}| + m\cdot\dim \langle R^{(\theta)}\rangle\nonumber\\
&\geq & 1 + m\cdot \dim \langle R^{(\theta)}\rangle.\nonumber
\end{eqnarray}
But this is a contradiction to $T$ being a qualified $mk$-system, i.e.
to \eqref{2.4}. Therefore $N\in\Z_{>0}$.

It rests to construct a new decomposition $T=S^{(1)}+...+S^{(m)}$ with 
\eqref{2.6}. For this,  we show now that a sequence $(b_j)_{j=0,...,N}$
of elements of $\supp T$ with the following properties can be chosen:
\begin{eqnarray}\label{2.12}
b_j&\in &R^{(j)}\qquad\textup{for }0\leq j\leq N,\\
b_j&\notin& R^{(j-m)}\qquad\textup{for }m+1\leq j\leq N,\label{2.13}\\
\lambda_{[j]}(b_{j-1},b_j)&\neq& 0 \qquad\textup{for }1\leq j\leq N,\label{2.14}\\
b_N&\notin& \langle T^{([N+1])}\rangle .\label{2.15}
\end{eqnarray}
We construct the elements in the order $b_N,b_{N-1},...,b_0$ and
start with $b_N$: By construction $R^{(N)}\not\subset 
\langle T^{([N+1])}\rangle$. Choose $b_N\in R^{(N)}-\langle T^{([N+1])}\rangle$.
\eqref{2.12} and \eqref{2.15} hold.
If $N\geq m+1$ then $R^{(N-m)}\subset \langle T^{([N+1])}\rangle$,
so \eqref{2.13} holds.

If $b_N,...,b_j$ are constructed with \eqref{2.12} -- \eqref{2.15} 
for some $j\geq 1$, then choose $b_{j-1}\in R^{(j-1)}$ with \eqref{2.14}.
This is possible because of \eqref{2.9}. If $b_{j-1}\in R^{(j-1-m)}$
in the case $j-1\geq m+1$, then by \eqref{2.9} also $b_j\in R^{(j-m)}$,
a contradiction. This shows the existence of $b_N,b_{N-1},...,b_0$ as claimed.
They are all different because of \eqref{2.12} and \eqref{2.13}.

Now define a sequence $(\www R^{(j)})_{j=1,...,N}$ of subsets of $\supp T$
as follows.
\begin{eqnarray}\label{2.16}
\www R^{(j)}:=\Bigl(R^{(j)}-\bigcup_{i\in\{1,...,j\}:[i]=j}\{b_i\}\Bigr) 
\cup \bigcup_{i\in\{0,...,j-1\}:[i+1]=j}\{b_i\},
\end{eqnarray}
so that
\begin{eqnarray}\label{2.17}
\www R^{(j)}&=& \Bigl(\bigl((R^{(j)}-\{b_j\})\cup\{b_{j-1}\}\bigr) 
-R^{(j-m)}\Bigr)\cup \www R^{(j-m)}\\
&&\textup{for }j\geq m+1.\nonumber
\end{eqnarray}
We claim that for all $j\in\{1,...,N\}$
\begin{eqnarray}\label{2.18}
\langle \www R^{(j)}\rangle = \langle R^{(j)}\rangle
\end{eqnarray}
holds. For $1\leq j\leq \min(m,N)$ this follows from \eqref{2.14}.
If $N\geq m+1$, then for $m+1\leq j\leq N$ this follows inductively: 
The induction hypothesis gives
\begin{eqnarray}\label{2.19}
\langle \www R^{(j-m)}\rangle = \langle R^{(j-m)}\rangle
\end{eqnarray}
Thus
\begin{eqnarray}
\langle \www R^{(j)}\rangle &=&
\left\langle\Bigl(\bigl((R^{(j)}-\{b_j\})\cup\{b_{j-1}\}\bigr) 
-R^{(j-m)}\Bigr)\cup \www R^{(j-m)}\right\rangle  \nonumber\\
&=& \left\langle\Bigl(\bigl((R^{(j)}-\{b_j\})\cup\{b_{j-1}\}\bigr) 
-R^{(j-m)}\Bigr)\right\rangle + \langle\www R^{(j-m)}\rangle \nonumber\\
&=& \left\langle\Bigl(\bigl((R^{(j)}-\{b_j\})\cup\{b_{j-1}\}\bigr) 
-R^{(j-m)}\Bigr)\right\rangle + \langle R^{(j-m)}\rangle \nonumber\\
&=& \left\langle\bigl((R^{(j)}-\{b_j\})\cup\{b_{j-1}\}\bigr) 
\right\rangle \nonumber\\
&=& \langle R^{(j)}\rangle .\nonumber
\end{eqnarray}
The last equality uses \eqref{2.14}.

Because $R^{(N)}\not\subset \langle T^{([N+1])}\rangle$,
$\langle T^{([N+1])}\rangle\subsetneqq V$ and 
$\supp T^{([N+1])}-B^{([N+1])}\neq \emptyset$. Choose
\begin{eqnarray}\label{2.20}
b_{N+1}\in\supp T^{([N+1])}-B^{([N+1])}
\end{eqnarray}
arbitrarily

Now we will define $k$-systems $S^{(1)},...,S^{(m)}$ with
\begin{eqnarray}\label{2.21}
S^{(1)}+...+S^{(m)}=T
\end{eqnarray}
and \eqref{2.6}. We have to distinguish two cases.

{\bf First case, $N\leq m-1$:}
\begin{eqnarray}\label{2.22}
S^{(j)}&:=& T^{(j)}-[b_j]+[b_{j-1}]
\qquad\textup{for }2\leq j\leq N+1,\\
S^{(1)}&:=& T^{(1)}-[b_1]+[b_{N+1}],\label{2.23}\\
S^{(j)}&:=& T^{(j)}\qquad\textup{for }N+2\leq j\leq m.\label{2.24}
\end{eqnarray}
These are $m$ $k$-systems.
\eqref{2.21} holds obviously. The subspaces generated by
the $k$-systems $S^{(j)}$ are as follows.
\begin{eqnarray*}
\langle S^{(j)}\rangle &:=& \langle T^{(j)}\rangle 
\qquad\textup{for }2\leq j\leq N,\quad \textup{by \eqref{2.14}},\\
\langle S^{(N+1)}\rangle &:=& \langle T^{(N+1)}\rangle +\langle b_N\rangle \supsetneqq \langle T^{(N+1)}\rangle\quad\textup{by \eqref{2.20} and \eqref{2.15}},\\
\langle S^{(1)}\rangle &:=& \langle T^{(1)}\rangle +
 \langle b_{N+1}\rangle\supset  \langle T^{(1)}\rangle 
 \quad \textup{by \eqref{2.14}},\\
\langle S^{(j)}\rangle &:=& \langle T^{(j)}\rangle 
\qquad\textup{for }N+2\leq j\leq m.
\end{eqnarray*}
Together these give \eqref{2.6}.

{\bf Second case, $N\geq m$:}
Define
\begin{eqnarray}\label{2.25}
S^{(j)}&:=& T^{(j)} -\sum_{i\in\{1,..,N+1\}:\ [i]=j}[b_i]
+\sum_{i\in\{1,..,N\}:\ [i+1]=j}[b_i]\quad\text{for }
2\leq j\leq m,\hspace*{1cm}\\
S^{(1)} &:=& T^{(1)} -\sum_{i\in\{1,..,N+1\}:\ [i]=1}[b_i]
+\sum_{i\in\{1,..,N\}:\ [i+1]=1}[b_i]+[b_{N+1}].\label{2.26}
\end{eqnarray}
These are $m$ $k$-systems.
\eqref{2.21} holds because of \eqref{2.25} and \eqref{2.26}.
For $j\in\{1,...,m\}$ let $\alpha(j)$ be the unique integer
in $\{N-m+1,N-m+2,...,N\}$ with $[\alpha(j)]=j$.
Especially $\alpha([N+1])=N-m+1$.
Then for $j\in\{1,...,m\}-\{1,[N+1]\}$
\begin{eqnarray}\label{2.27}
S^{(j)}&=& T^{(j)} -\sum_{b\in R^{(\alpha(j))}}[b]
+\sum_{b\in\www R^{(\alpha(j))}}[b]
\end{eqnarray}
and
\begin{eqnarray}\label{2.28}
\langle S^{(j)}\rangle &=& \langle T^{(j)}\rangle 
\quad \textup{by \eqref{2.18}}.
\end{eqnarray}
If $1\neq [N+1]$ then
\begin{eqnarray}\label{2.29}
S^{(1)} &=& T^{(1)} -\sum_{b\in R^{(\alpha(1))}}[b]
+(\sum_{b\in\www R^{(\alpha(1))}}[b]-[b_0])+[b_{N+1}],\\
S^{([N+1])} &=& T^{([N+1])} -\sum_{b\in R^{(N-m+1)}}[b]
+\sum_{b\in\www R^{(N-m+1)}}[b]-[b_{N+1}]+[b_N]\label{2.30}
\end{eqnarray}
and 
\begin{eqnarray}
\langle S^{(1)}\rangle &=& \langle T^{(1)}\rangle +
 \langle b_{N+1}\rangle\supset  \langle T^{(1)}\rangle 
\quad \textup{by \eqref{2.18}},\label{2.31}\\
\langle S^{([N+1])}\rangle &=& \langle T^{([N+1])}\rangle +\langle b_N\rangle 
\supsetneqq \langle T^{([N+1])}\rangle\label{2.32}\\
&&\qquad\textup{by \eqref{2.18}, \eqref{2.20} and \eqref{2.15}}.\nonumber
\end{eqnarray}
If $1=[N+1]$ then
\begin{eqnarray}\label{2.33}
S^{(1)} &=& T^{(1)} -\sum_{b\in R^{(N-m+1)}}[b]
+(\sum_{b\in\www R^{(N-m+1)}}[b]-[b_0])+[b_N]
\end{eqnarray}
(in the case $1=[N+1]$, $b_{N+1}$ is not used) and
\begin{eqnarray}\label{2.34}
\langle S^{(1)}\rangle &=& \langle T^{(1)}\rangle +
 \langle b_N\rangle\supsetneqq  \langle T^{(1)}\rangle 
\quad\textup{by \eqref{2.18} and \eqref{2.15}}.
\end{eqnarray}
In both cases, $1\neq [N+1]$ and $1=[N+1]$, \eqref{2.6} holds.
This finishes the first step and proves part (a) of theorem
\ref{t2.4} for $l=0$.

\medskip
{\bf Second step:}
The cases $l\geq 1$ are proved by induction in $l$.
Fix some $l\geq 1$. Suppose that any qualified 
$(m\www k+\www l)$-system
with $\www k\in\Z_{>0}$ arbitrary and 
$\www l\in\{0,1,...,l-1\}$ is strong.
Let $T$ be a qualified $(mk+l)$-system.
As in the beginning of the proof we can suppose
$0\notin\supp T$ and $T(b)=1$ for all $b\in\supp T$.

Define $g:=\dim\langle A_2(T)\rangle$.
Choose an arbitrary element $a\in A_2(T)$.
The subsystem
\begin{eqnarray}\label{2.35}
R:=\sum_{b\in A_2(T)-\{a\}}[b]
\end{eqnarray}
of $T$ is an $(mg+l-1)$-system
by definition of $A_2(T)$. Furthermore, it is a qualified
$(mg+l-1)$-system with respect to the vector space 
$\langle A_2(T)\rangle$, again by
definition of $A_2(T)$: 
For any $U\subsetneqq \langle A_2(T)\rangle$, 
$\mu(T,U)\leq l-1+m\dim U$ by the minimality of 
$\langle A_2(T)\rangle$, thus also $\mu(R,U)\leq l-1+m\dim U$. 
For $U=\langle A_2(T)\rangle$,
this holds because $R$ is an $(mg+l-1)$-system.

By induction hypothesis, $R$ is a strong
$(mg+l-1)$-system with respect to the vector space 
$\langle A_2(T)\rangle$, so it has a
strong decomposition $R=R^{(1)}+...+R^{(m+1)}$.
Therefore $R+[a]=\sum_{b\in A_2(T)}[b]$ is a strong
$(mg+l)$-system with the strong decomposition
$R+[a]=R^{(1)}+...+R^{(m)}+(R^{(m+1)}+[a])$.

Consider the system 
\begin{eqnarray}\label{2.36}
S:=\sum_{b\in \supp T-A_2(T)}[b+\langle A_2(T)\rangle]
\end{eqnarray}
of vectors in the quotient space $V/\langle A_2(T)\rangle$.
It is a qualified $m(k-g)$-system 
with respect to the vector space 
$V/\langle A_2(T)\rangle$
because $T$ is a qualified $(mk+l)$-system and 
$R+[a]$ is a qualified $(mg+l)$-system.
By the first step in this proof, $S$ has a strong decomposition
$S=S^{(1)}+...+S^{(m)}$.

The vectors in $\supp S\subset V/\langle A_2(T)\rangle$ 
lift uniquely to vectors in $\supp T$, because
$T(b)=1$ for all $b\in\supp T$. Let 
$\www S=\www S^{(1)}+...+\www S^{(m)}$ be the corresponding lift
to $V$ of $S$ and its decomposition.
Then
\begin{eqnarray}\label{2.37}
T=(R^{(1)}+\www S^{(1)})+...+(R^{(m)}+\www S^{(m)})
+(R^{(m+1)}+[a])
\end{eqnarray}
is a strong decomposition of $T$. This finishes the second step
and the proof of part (a).

\medskip
(b) Because of lemma \ref{t2.3} (c), it rests to show
$A_2(T)\subset A_1(T)$. But this follows from the 
second step above and especially the strong decomposition
\eqref{2.37} of $T$. There $a\in A_2(T)$ was arbitrary.
\hfill$\Box$

\begin{remarks}\label{t2.5}
Part (a) of theorem \ref{t2.4} has some
similarity with the marriage theorem of Hall:

{\it Let $A$ and $B$ be nonempty finite sets with
$|A|\leq |B|$, and let $f:A\to\PP(B)(:=$ the set of subsets
of $B$) be a map. 
Then a map $g:A\to B$ with $g(a)\in f(a)$ exists if
and only 
\begin{eqnarray}\label{2.38}
|\bigcup_{c\in C}f(c)|\geq |C|\qquad \textup{for all }
C\subset A\textup{ with }C\neq\emptyset.
\end{eqnarray}
}

In theorem \ref{t2.4} and in the marriage theorem of Hall,
the conditions \eqref{2.4} respectively \eqref{2.38}
are obviously necessary, but that they are sufficient
requires a nontrivial proof.
\end{remarks}

\section{An equivalence between index systems}\label{s3}
\setcounter{equation}{0}

\noindent
This section prepares the proof of the main result, theorem \ref{t4.2}.
It builds on section \ref{s2}.

Start with three positive integers $k$ and $n$ and $m$ with $k<n$
and $m\geq 2$, 
with a field $K$, a $K$-vector space $V$ of dimension $k$ and
a map
\begin{eqnarray}\label{3.1}
v:J\to V,\quad i\mapsto v(i)=:v_i,\quad\textup{ where }J:=\{1,...,n\},
\end{eqnarray}
with the property
\begin{eqnarray}\label{3.2}
\langle v_1,...,v_n\rangle =V.
\end{eqnarray}

\begin{definition}\label{t3.1}
\begin{list}{}{}
\item[(a)] A system $T=\sum_{i\in J}T(i)\cdot[i]$ of elements of $J$ induces
the system
\begin{eqnarray}\label{3.3}
v^{sys}(T):=\sum_{i\in J}T(i)\cdot [v_i]
\end{eqnarray}
of elements of $V$. Of course $|T|=|v(T)|$.

\item[(b)]
A {\it strong decomposition} of an $(mk+l)$-system $T$ 
of elements of $J$ for $l\in\{0;1\}$
is a decomposition $T=T^{(1)}+...+T^{(m+1)}$ into $m+1$ systems
such that the induced decomposition $v^{sys}(T)=v^{sys}(T^{(1)})
+...+v^{sys}(T^{(m+1)})$ of $v^{sys}(T)$ is strong
(definition \ref{t2.2} (c)). 

\item[(c)] 
An $(mk+l)$-system $T$ of elements of $J$ for $l\in\{0;1\}$
is {\it strong} if it admits a strong decomposition.
Of course, this holds if and only if the system $v^{sys}(T)$ 
is strong (definition \ref{t2.2} (d)).

\item[(d)]
A {\it good decomposition} of an $N$-system $T$ of elements of $J$
with $N\geq mk+1$ is a decomposition $T=T_1+T_2$ 
into two systems such that $T_2$ is a strong $(mk+1)$-system 
of elements of $J$.

\item[(e)]
Two good decompositions $T_1+T_2=T$ and $S_1+S_2=T$ 
of an $N$-system $T$ of elements of $J$ with $N\geq mk+1$
are {\it locally related}, 
notation: $(S_1,S_2)\sim_{loc} (T_1,T_2)$, 
if there are strong decompositions 
$S^{(1)}_2+...+S^{(m+1)}_2=S_2$ of $S_2$ and 
$T^{(1)}_2+...+T^{(m+1)}_2=T_2$ of $T_2$ with 
$S^{(j)}_2=T^{(j)}_2$ for $1\leq j\leq m$. 
Of course, $\sim_{loc}$ is a reflexive and symmetric relation.

\item[(f)]
Two good decompositions $T_1+T_2=T$ and $S_1+S_2=T$ 
of an $N$-system $T$ of elements of $J$ with $N\geq mk+1$
are {\it equivalent}, 
notation: $(S_1,S_2)\sim (T_1,T_2)$, 
if there is a sequence $\sigma_1,\sigma_2,...,\sigma_r$
for some $r\in\Z_{\geq 1}$ of good decompositions of $T$
such that $\sigma_1=(S_1,S_2)$, $\sigma_r=(T_1,T_2)$ and
$\sigma_j\sim_{loc}\sigma_{j+1}$ for $j=1,...,r-1$.
Of course, $\sim$ is an equivalence relation.

\item[(g)]
The distance $d_H(S,T)$ between two 
systems $S$ and $T$ of elements of $J$ is the number
\begin{eqnarray}\label{3.4}
d_H(S,T):=\sum_{i\in J}|S(i)-T(i)|\in\Z_{\geq 0}.
\end{eqnarray}
This defines a metric on the set of systems of elements of $J$.
\end{list}
\end{definition}

The main result of this section is the following theorem \ref{t3.2}.

\begin{theorem}\label{t3.2}
Let $T$ be an $N$-system of elements of $J$ for some $N\geq mk+1$
which has good decompositions. Then all its good decompositions
are equivalent.
\end{theorem}

The theorem will be proved after the proof of lemma \ref{t3.3}.

\begin{lemma}\label{t3.3}
Let $S$ and $T$ be two strong $(mk+1)$-systems of elements of $J$.
At least one of the following two alternatives holds.

\begin{list}{}{}
\item[$(\alpha)$] 
$T$ has a strong decomposition $T=T^{(1)}+...+T^{(m+1)}$
with $T^{(m+1)}=[i]$ for some $i\in \supp T$ with $T(i)>S(i)$.
\item[$(\beta)$] 
$T$ and $S$ have strong decompositions
$T=T^{(1)}+...+T^{(m+1)}$ and $S=S^{(1)}+...+S^{(m+1)}$
with $T^{(m+1)}=S^{(m+1)}$.
\end{list}
\end{lemma}

The lemma builds on section \ref{s2}, especially on part (b)
of theorem \ref{t2.4}.

{\bf Proof of lemma \ref{t3.3}:}
Define $A_1(T):=\{i\in \supp T\, |\, v_i\in A_1(v^{sys}(T))\}$
and analogously $A_1(S)$. Then
\begin{eqnarray}
A_1(T)&=&\{i\in\supp T\, |\, \exists\ 
\textup{ a strong decomposition }
T^{(1)}+...+T^{(m+1)}\nonumber \\
&& \textup{ with }T^{(m+1)}=[i]\},\label{3.5}\\
\sum_{i\in A_1(T)}T(i)&=&1+m\dim \langle A_1(v^{sys}(T))\rangle \label{3.6}\\
&=&1+m\dim\langle v_i\, |\, i\in A_1(T)\rangle, \nonumber
\end{eqnarray}
and analogously for $A_1(S)$. Here \eqref{3.6} follows
from theorem \ref{t2.4} (b).

Suppose that $(\alpha)$ does not hold.
Then for any $i\in A_1(T)$ $S(i)\geq T(i)>0$. Especially
\begin{eqnarray}\label{3.7}
\sum_{i\in A_1(T)}S(i)\geq 1+m\dim \langle v_i\, |\, i\in A_1(T)\rangle.
\end{eqnarray}
\eqref{3.7} and the argument in the proof of lemma
\ref{t2.3} (c) show $A_1(S)\subset A_1(T)$
(and $S(i)=T(i)$ for $i\in A_1(T)$). Thus $(\beta)$ holds.
\hfill$\Box$

\bigskip
{\bf Proof of theorem \ref{t3.2}:}
Let $(S_1,S_2)$ and $(T_1,T_2)$ be two different good decompositions
of an $N$-system $T$ of elements of $J$ (with $N\geq mk+1$).
Then $S_2$ and $T_2$ are strong $(mk+1)$-systems of elements of $J$.
At least one of the two alternatives $(\alpha)$ and $(\beta)$ in
lemma \ref{t3.3} holds for $S_2$ and $T_2$.

\medskip
{\bf First case, $(\alpha)$ holds:} 
Let $T_2=T_2^{(1)}+...+T_2^{(m+1)}$ be a strong decomposition
with $T_2^{(m+1)}=[i]$ for some $i\in\supp T_2$ with 
$T_2(i)>S_2(i)$.
Then a $j\in \supp T$ with $T_1(j)>S_1(j)$ and $T_2(j)<S_2(j)$ exists. 
The decomposition 
\begin{eqnarray}\label{3.8}
T=R_1+R_2\quad\textup{with }R_1=T_1-[j]+[i],\quad R_2=T_2+[j]-[i]
\end{eqnarray}
is a good decomposition of $T$ because 
$T_2^{(1)}+...+T_2^{(m)}+[j]$ is a strong decomposition of $R_2$. 
The good decompositions $(R_1,R_2)$ and $(T_1,T_2)$ are locally
related, $(R_1,R_2)\sim_{loc}(T_1,T_2)$, and thus equivalent,
\begin{eqnarray}\label{3.9}
(R_1,R_2)\sim(T_1,T_2).
\end{eqnarray}
Furthermore, 
\begin{eqnarray}\label{3.10}
d_H(R_2,S_2)=d_H(T_2,S_2)-2.
\end{eqnarray}

\medskip
{\bf Second case, $(\beta)$ holds:}
Let $T_2=T_2^{(1)}+...+T_2^{(m+1)}$ and 
$S_2=S_2^{(1)}+...+S_2^{(m+1)}$
be strong decompositions of $T_2$ and $S_2$ with 
$T_2^{(m+1)}=S_2^{(m+1)}=[a]$ for some $a\in \supp T$.
Two elements $b,c\in\supp T$ with 
$T_1(b)>S_1(b)$, $T_2(b)<S_2(b)$, and $T_1(c)<S_1(c)$, $T_2(c)>S_2(c)$
exist. 
Consider the decompositions of $T$ and $S$,
\begin{eqnarray}\label{3.11}
T&=&R_1+R_2\quad\textup{with }R_1=T_1-[b]+[a],R_2=T_2+[b]-[a],\\
S&=&Q_1+Q_2\quad\textup{with }Q_1=S_1-[c]+[a],Q_2=S_2+[c]-[a].\label{3.12}
\end{eqnarray}
They are good decompositions because $R_2$ has the strong
decomposition $R_2=T^{(1)}+...+T^{(m)}+[b]$ 
and $Q_2$ has the strong decomposition
$Q_2=S^{(1)}+...+S^{(m)}+[c]$. 
The local relations
\begin{eqnarray*}
(R_1,R_2)\sim_{loc} (T_1,T_2)\quad\textup{and}\quad 
(Q_1,Q_2)\sim_{loc}(S_1,S_2)
\end{eqnarray*}
and the equivalences
\begin{eqnarray}\label{3.13}
(R_1,R_2)\sim (T_1,T_2)\quad\textup{and}\quad 
(Q_1,Q_2)\sim (S_1,S_2)
\end{eqnarray}
hold. Furthermore
\begin{eqnarray}\label{3.14}
d_H(R_2,Q_2)=d_H(T_2,S_2)-2.
\end{eqnarray}

\medskip
The properties \eqref{3.10}, \eqref{3.11}, \eqref{3.13} and \eqref{3.14}
show that in both cases the equivalence classes of $(S_1,S_2)$ and
$(T_1,T_2)$ contain good decompositions whose second members
are closer to one another with respect to the metric $d_H$
than $T_2$ and $S_2$. This shows that $(S_1,S_2)$ and $(T_1,T_2)$
are in one equivalence class.  \hfill$\Box$

\section{Potentials of the first and second kind}\label{s4}
\setcounter{equation}{0}

The main part of this section is devoted to the proof of
theorem \ref{t1.2}. At the end some remarks on the
relation to families of arrangements, 
Frobenius manifolds, F-manifolds and possible
extensions of the work here are made.

\begin{remark}\label{t4.1}
This is a reminder of the notion of a Higgs field and the 
meaning of the condition $\nnn^K(C)=0$ in definition \ref{t1.1}.
There a Higgs field is an $\OO_M$-linear map
\begin{eqnarray}\label{4.1}
C:\OO(K)\to \Omega^1_K\otimes \OO(K)
\end{eqnarray}
such that all the endomorphisms $C_X,X\in\TT_M$, commute.
$\nnn^K(C)\in \Omega^2_M\otimes\OO(K)$ is the 2-form on $M$
with values in $K$ such that for $X,Y\in\TT_M$
\begin{eqnarray}\label{4.2}
\nnn^K(C)(X,Y)&=& \nnn^K_X(C_Y)-\nnn^K_Y(C_X)-C_{[X,Y]}.
\end{eqnarray}
Now $\nnn^K(C)=0$ says
\begin{eqnarray}\label{4.3}
\nnn^K_{\paa_i}(C_{\paa_j})=\nnn^K_{\paa_j}(C_{\paa_i}).
\end{eqnarray}
\end{remark}

{\bf Proof of theorem \ref{t1.2}:}
Let $(M,K,\nnn^K,C,S,\zeta,V,(v_1,...,v_n))$ be a Frobenius 
like structure of some order $(n,k,m)\in\Z_{>0}^3$.

We need some notations. If $T\in\Z_{\geq 0}[J]$ is a system
of elements of $J$, then 
\begin{eqnarray*}
(z-x)^T&:=&\prod_{i\in J}(z_i-x_i)^{T(i)}\quad\textup{for any }x\in\C^n,\\
T!&:=&\prod_{i\in J}T(i)!,\\
\paa_T&:=&\prod_{i\in J}\paa_{z_i}^{T(i)},\\
C_T&:=&\prod_{i\in J}C_{\paa_{z_i}}^{T(i)}.
\end{eqnarray*}
Thus, if $S$ and $T$ are systems of elements of $J$, then
\begin{eqnarray}
\paa_T(z-x)^S=\left\{\begin{array}{ll}
0&\textup{ if }T\not\leq S,\\
\frac{S!}{(S-T)!}\cdot (z-x)^{S-T}& \textup{ if }T\leq S,
\end{array}\right. \label{4.4}
\end{eqnarray}
for any $x\in\C^n$.

\medskip
The existence of a (not just local, but even global) 
potential $Q$ of the first kind is trivial.
The function
\begin{eqnarray}\label{4.5}
Q&:=& \sum_{T\textup{ with }(*)}\frac{1}{T!}\cdot 
S(C_T \zeta,\zeta,...,\zeta)\cdot z^T \quad(m\textup{ times }\zeta),
\hspace*{1cm}\\
(*)&:& T\in\Z_{\geq 0}[J]
\textup{ is a strong }mk\textup{-system (definition \ref{t3.1}(c))}.\nonumber
\end{eqnarray}
works. It is a homogeneous polynomial of degree $mk$
and contains only monomials which are relevant for \eqref{1.2}.
In fact, one can add to this $Q$ an arbitrary linear combination
of the monomials $z^T$ for the $mk$-systems $T$ which are not strong,
so which are not relevant for \eqref{1.2}.

\medskip
The existence of a potential $L$ of the second kind is not trivial.
Let some $x\in M$ be given. We make the power series ansatz
\begin{eqnarray}\label{4.6}
L&:=& \sum_{T\in\Z_{\geq 0}[J]} a_T\cdot (z-x)^T,
\end{eqnarray}
where the coefficients $a_T$ have to be determined.
If $T$ satisfies $|T|\leq mk$ or if it satisfies
$|T|\geq mk+1$, but does not admit a good decomposition
(definition \ref{t3.1} (d)), then the conditions \eqref{1.3}
are empty for $a_T(z-x)^T$ because of \eqref{4.4},
so then $a_T$ can be chosen arbitrarily, e.g. $a_T:=0$ works.

Now consider $T$ with $|T|\geq mk+1$ which admits good decompositions.
Then each good decomposition $T=T_1+T_2$ gives 
via \eqref{1.3} a candidate
\begin{eqnarray}\label{4.7}
a_T(T_1,T_2)&:=& \frac{1}{T!}\cdot \left(\paa_{T_1}
S(C_{T_2}\zeta,\zeta,...,\zeta)\right)(x),
\end{eqnarray}
for the coefficient $a_T$ of $(z-x)^T$ in $L$.
We have to show that the candidates $a_T(T_1,T_2)$ for all
good decompositions $(T_1,T_2)$ of $T$ coincide.

Suppose that two good decompositions $(T_1,T_2)$ and $(S_1,S_2)$
are locally related, $(T_1,T_2)\sim_{loc}(S_1,S_2)$
(definition \ref{t3.1} (e)), but not equal.
Then there are strong decompositions $T_2=T_2^{(1)}+...+T_2^{(m)}+[a]$
and $S_2=T_2^{(1)}+...+T_2^{(m)}+[b]$ with $a\neq b$,
and thus also $T_1-[b]=S_1-[a]\in\Z_{\geq 0}[J]$ holds.
Because any $T_2^{(j)}$, $j\in\{1,...,m\}$, is independent,
$C_{T_2^{(j)}}\zeta$ is $\nnn^K$-flat.
This and \eqref{4.3} give
\begin{eqnarray}
&&\paa_{z_b}S(C_{T_2}\zeta,\zeta,...,\zeta)\nonumber\\
&=& \paa_{z_b}S(C_{\paa_{z_a}}C_{T_2^{(1)}}\zeta, C_{T_2^{(2)}}\zeta,...,
C_{T_2^{(m)}}\zeta)\nonumber\\
&=& S(\nnn^K_{\paa_{z_b}}(C_{\paa_{z_a}})C_{T_2^{(1)}}\zeta, C_{T_2^{(2)}}\zeta,...,
C_{T_2^{(m)}}\zeta)\nonumber\\
&=& S(\nnn^K_{\paa_{z_a}}(C_{\paa_{z_b}})C_{T_2^{(1)}}\zeta, C_{T_2^{(2)}}\zeta,...,
C_{T_2^{(m)}}\zeta)\nonumber\\
&=& \paa_{z_a}S(C_{\paa_{z_b}}C_{T_2^{(1)}}\zeta, C_{T_2^{(2)}}\zeta,...,
C_{T_2^{(m)}}\zeta)\nonumber\\
&=&\paa_{z_a}S(C_{S_2}\zeta,\zeta,...,\zeta).
\label{4.8}
\end{eqnarray}
This implies 
\begin{eqnarray}\label{4.9}
a_T(T_1,T_2) = a_T(S_1,S_2),
\end{eqnarray}
so the locally related good decompositions $(T_1,T_2)$ and $(S_1,S_2)$
give the same candidate for $a_T$. 
Thus all equivalent (definition \ref{t3.1} (f)) good decompositions
give the same candidate for $a_T$.
By theorem \ref{t3.2}, all good decompositions of $T$ are equivalent.
Therefore they all give the same candidate for $a_T$.
Thus a potential $L$ of the second kind exists as a formal power series
as in \eqref{4.6}.

It is in fact a convergent power series because of the following.
There are finitely many strong $mk$-systems $T_2$. Each determines
the coefficients $a_T$ for all $T\geq T_2$. 
We put $a_T:=0$ for $T$ which do not admit
good decompositions. The part of $L$ in \eqref{4.6} which is determined
by some strong $mk$-system $T_2$ is a convergent power series.
Thus $L$ is the {\it union} of finitely many overlapping
convergent power series. It is easy to see that it is itself convergent. 
This finishes the proof of theorem \ref{t1.2}.
\hfill$\Box$

\begin{remark}\label{t4.2}
In \cite[ch. 3]{V2} families or arrangements are considered which
give rise to Frobenius like structures 
$(M,K,\nnn^K,C,S,\zeta,V,(v_1,...,v_n))$ of order $(n,k,2)$,
see the special case of generic arrangements in \cite{V1,V3}.

Start with two positive integers $k$ and $n$ with $k<n$
and with a matrix $B:=(b_i^j)_{i=1,..,n;j=1,..,k}\in M(n\times k,\C)$
with $\rank B=k$. Define $J:=\{1,...,n\}$. 
Here the vector space $V$ and the vectors $v_1,...,v_n$ are
\begin{eqnarray}\label{4.10}
V&=&M(1\times k,\C),\\
v_i&=& (b^j_i)_{j=1,...k}\in V\qquad\textup{for }i=1,...,n.
\end{eqnarray}
We assume that $B$ is such that all vectors $v_i$ are nonzero.

Consider $\C^n\times \C^k$ with the
coordinates $(z,t)=(z_1,...,z_n,t_1,...,t_k)$ and with the projection
$\pi:\C^n\times \C^k\to\C^n$. Define the functions
\begin{eqnarray}\label{4.11}
g_i:=\sum_{j=1}^kb_i^j\cdot t_j,\quad f_i:=g_i+z_i
\quad\textup{for }i\in J
\end{eqnarray}
on $\C^n\times \C^k$. 

We obtain on $\C^n\times \C^k$ the arrangement 
$\CC=\{H_i\}_{i\in J}$, where $H_i$ is the zero set of $f_i$.
Let $U(\CC):=\C^n\times \C^k-\bigcup_{i\in J}H_i$ be the complement.
For every $x\in\C^n$, the arrangement $\CC$ restricts to an 
arrangement $\CC(x)$ on $\pi^{-1}(x)\cong\C^k$. 
For almost all $x\in\C^k$ the arrangement $\CC(x)$ is {\it essential}
(definition in \cite{V2})
with normal crossings. The subset $\Delta\subset\C^n$ where
this does not hold, is a hypersurface and is called the {\it discriminant},
see \cite[3.2]{V2}. Define $M:=\C^n-\Delta$.

A set $I=\{i_1,...,i_k\}\subset J$ is independent, i.e.
$(v_{i_1},...,v_{i_k})$ is a basis of $V$, if and only if
(for some or equivalently for any $x\in \C^n$) 
the hyperplanes $H_{i_1}(x),...,H_{i_k}(x)$ are transversal.

Let $a=(a_1,...,a_n)\in(\C^*)^n$ be a system of {\it weights}
such that for any $x\in M$ the weighted arrangement
$(\CC(x),a)$ is {\it unbalanced}: See \cite{V2} for the definition
of {\it unbalanced}, e.g. $a\in\R_{>0}^n$ is unbalanced,
also a generic system of weights is     unbalanced.
The {\it master function} of the weighted arrangement $(\CC,a)$
is
\begin{eqnarray}\label{4.12}
\Phi_a(z,t):=\sum_{i\in J}a_i\log f_i.
\end{eqnarray}
Several deep facts are related to this master function.
We use some of them in the following. See \cite{V2} for references.

For $z\in M$ all critical points of $\Phi_{a}$ are isolated,
and the sum $\mu$ of their Milnor numbers is independent of 
the unbalanced weight $a$ and the parameter $z\in M$. 
The bundle 
\begin{eqnarray}\label{4.13}
K:=\bigcup_{z\in M}K_z\quad\textup{with }
K_z:=\OO(U(\CC)\cap\pi^{-1}(z))/\left(\frac{\paa\Phi_a}{\paa t_j}\, |\,
j=1,...,k\right)
\end{eqnarray}
over $M$ is a vector bundle of $\mu$-dimensional algebras.

It comes equipped with the section $\zeta$ of unit elements 
$\zeta(z)\in K_z$, a Higgs field $C$, 
a {\it combinatorial connection} $\nnn^K$ and a pairing $S$.
The Higgs field $C:\OO(K)\to \Omega^1_M\otimes \OO(K)$
is defined with the help of the period map
\begin{eqnarray}\label{4.15}
\Psi:TM\to K,\quad \paa_{z_i}\mapsto \left[\frac{\paa\Phi_a}{\paa z_i}\right]
=\left[\frac{a_i}{f_i}\right]=:p_i
\end{eqnarray}
by
\begin{eqnarray}\label{4.16}
C_{\paa_{z_i}}(h):=p_i\cdot h\qquad\textup{ for }h\in K_z.
\end{eqnarray}
Because of 
\begin{eqnarray}\label{4.17}
0=\left[\frac{\paa\Phi_a}{\paa t_j}\right]
=\sum_{i=1}^n b^j_i p_i,
\end{eqnarray}
the Higgs field vanishes on the vector fields
$X_j:=\sum_{i=1}^n b^j_i\paa_i$, $j\in\{1,...,k\}$, 
\begin{eqnarray}\label{4.18}
C_{X_j}=0\qquad\textup{for }j\in\{1,...,k\}.
\end{eqnarray}
In fact the whole geometry of the family of arrangements is 
invariant with respect to the flows of these vector fields.

The sections $\det(b_i^j)_{i\in I,j=1,...,k}\cdot C_I\zeta$ 
for all independent sets $I=\{i_1,...,i_k\}\subset J$
generate the bundle $K$, and they satisfy only relations with 
constant coefficients in $\Z$. The combinatorial connection $\nnn^K$ is the unique
flat connection such that the sections $C_I\zeta$ for $I\subset J$
independent are $\nnn^K$-flat.
The sections $\det(b_i^j)_{i\in I,j=1,...,k}\cdot C_I\zeta$ for $I\subset J$
independent generate a $\nnn^K$-flat $\Z$-lattice structure on $K$.

The pairing $S$ comes from the Grothendieck residue with respect
to the volume form 
\begin{eqnarray}\label{4.19}
\frac{dt_1\land...\land dt_k}{\prod_{j=1}^k \frac{\paa\Phi_a}{\paa t_j}}.
\end{eqnarray}
It is a symmetric, nondegenerate, $\nnn^K$-flat, multiplication invariant
and Higgs field invariant.

\smallskip
The existence of potentials of the first and second kind for families of arrangements was conjectured in \cite{V1}. If all the $k\times k$ minors of the matrix $B=(b_i^j)$ are nonzero, the potentials were constructed in  \cite{V1}, cf. \cite{V3}. The potentials are given by explicit formulas in terms of the linear functions defining the hyperplanes in $\C^n$ composing the discriminant.

\end{remark}

\begin{remarks}\label{t4.3}
(i) The situation in remark \ref{t4.2} is in several aspects
richer than a Frobenius like structure of type $(n,k,m)$.
The bundle $K$ is a bundle of algebras.
The sections $C_I\zeta$ for independent sets $I\subset J$ generate
the bundle. The sections 
$\det(b_i^j)_{i\in I,j=1,...,k}\cdot C_I\zeta$
generate a flat $\Z$-lattice structure in $K$.
The Higgs field vanishes on the vector fields $X_1,...,X_k$.
The $m$-linear form $S$ is a pairing ($m=2$) and is nondegenerate.
We will not discuss the $\Z$-lattice structure, but we will discuss
some logical relations between the other enrichments and some 
implications of them.

\medskip
(ii) Let $(M,K,\nnn^K,C,S,\zeta,V,(v_1,...,v_n))$ be a Frobenius like
structure of order $(n,k,m)$. Suppose that it satisfies the 
{\it generation condition}
\begin{eqnarray}\label{4.20}
\text{(GC)}&& \textup{The sections }C_I\zeta\textup{ for independent sets }
I\subset J\\
&&\textup{generate the bundle }K. \nonumber
\end{eqnarray}
Let $\mu$ be the rank of $K$. Then for any $x\in M$, the endomorphisms
$C_X,X\in T_xM$, generate a $\mu$-dimensional commutative subalgebra
$A_z\subset\textup{End}(K_x)$. And any endomorphism which commutes
with them is contained in this subalgebra. This gives a rank $\mu$
bundle $A$ of commutative algebras. And the map
\begin{eqnarray}\label{4.21}
A\to K,\quad B\mapsto B\zeta,
\end{eqnarray}
is an isomorphism of vector bundles and induces a commutative and
associative multiplication on $K_x$ for any $x\in M$, with unit field
$\zeta(x)$. Therefore the special section $\zeta$  and 
the generation condition (GC), which exist and hold in remark 
\ref{t4.2}, give the multiplication on the bundle $K$ there.

\medskip
(iii) In the situation in (ii) with the condition (GC),
the $m$-linear form is multiplication invariant because it is 
Higgs field invariant. The condition (GC) implies also that it is symmetric:
\begin{eqnarray*}
S(C_{I_1}\zeta,C_{I_2}\zeta,...,C_{I_m}\zeta)
=S(C_{I_{\sigma(1)}}\zeta,C_{I_{\sigma(2)}}\zeta,...,C_{I_{\sigma(m)}}\zeta)
\end{eqnarray*}
for any independent sets $I_1,...,I_m$ and any permutation $\sigma\in S_m$.

\medskip
(iv)  For a Frobenius like structure $(M,K,\nnn^K,C,S,\zeta,V,(v_1,...,v_n))$ 
of order $(n,k,m)$ define the following $k$-dimensional space of linear combinations
of the coordinate vector fields $\paa_1,...,\paa_n$,
\begin{eqnarray}\label{4.22}
F^{inv}:=\{\sum_{i=1}^n\lambda(v_i)\paa_i\, |\, 
\lambda\in V^*\}.
\end{eqnarray}
It embeds into the tangent space $T_xM$ for any $x\in M$.
In remark \ref{t4.2} the space of these vector fields is the space
$\sum_{i=1}^k\C\cdot X_i$.

The {\it weak injectivity condition} is the condition
for any $x\in M$:
\begin{eqnarray}\label{4.23}
\text{(wIC)}&& \{X\in T_xM\, |\, C_X=0\}=F^{inv}\subset T_xM.
\end{eqnarray}
We expect that the potentials $Q$ and (locally) $L$ 
in theorem \ref{t4.2} can be chosen to be invariant 
with respect to the flows of the vector fields $X\in F^{inv}$ 
if the conditions (wIC) and (GC) hold. 
This is relevant for part (v) below.

\medskip
(v) Let $(M,K,\nnn^K,C,S,\zeta,V,(v_1,...,v_n))$ be a Frobenius like
structure of type $(n,k,2)$ with nondegenerate pairing $S$ which
satisfies the conditions (GC) and (wIC).
Consider an affine linear submanifold $N\subset M$ of dimension 
$n-k$ which is transversal to the orbits of $F^{inv}$.
One can identify it with the manifold of these orbits.
If potentials $Q$ and $L$ can be chosen to be constant on these orbits,
they live on this manifold.
An unfolding result in \cite{HM} can be applied to the restriction
of $(M,K,\nnn^K,C)$ to $N$. This will be discussed in the remarks \ref{t4.4}.
\end{remarks}

\begin{remarks}\label{t4.4}
(i) Let $(M,K,\nnn^K,C,S)$ be as follows.
$K\to M$ is a holomorphic vector bundle,
$\nnn^K$ is a holomorphic flat connection on $K$,
$C$ is a Higgs field on $K$ with $\nnn^K(C)=0$,
and $S$ is a holomorphic symmetric nondegenerate
$\nnn^K$-flat and Higgs field invariant pairing on $K$.
Let $p:\P^1\times M\to M$ be the projection.
The holomorphic vector bundle $H:=p^*K$ on $\P^1\times M$
is a family of trivial vector bundles $H|_{\P^1\times\{z\}},z\in M,$
on $\P^1$. Extend $\nnn^K$, $C$ and $S$ canonically to $H$.
Define 
\begin{eqnarray}\label{4.24}
\nnn&:=& \nnn^K+\frac{1}{\kappa}C,\\
\nnn&:&\OO(H)\to \OO_{\P^1\times M}(\{0\}\times M)\cdot\Omega_M^1\otimes
\OO(H),\nonumber
\end{eqnarray}
here $\kappa$ is the coordinate on $\C\subset\P^1$.
Then $\nnn$ restricts for any $\kappa\in\P^1-\{0\}$ to a flat connection
on $H|_{\{\kappa\}\times M}$ and has a pole of order 1 along 
$\{0\}\times M$.
Define a pairing $P:H_{(\kappa,z)}\times H_{(-\kappa,z)}\to\C$
for any $(\kappa,z)\in\P^1\times M$ by
\begin{eqnarray}\label{4.25}
P(a(\kappa,z),b(-\kappa,z))&:=& S(a(z),b(z)),
\end{eqnarray}
here $a(\kappa,z)\in H_{\kappa,z}$ and $b(-\kappa,z)\in H_{-\kappa,z}$
are the canonical lifts of elements $a(z),b(z)\in K_z$.
Then $P$ is a holomorphic symmetric nondegenerate $\nnn$-flat pairing.

In the notation of \cite{HM}, the tuple $(H\to\P^1\times M,\nnn,P)$
is a $(trTLP(0))$-structure. One can recover 
$(M,K,\nnn^K,C,S)$ from it. So one has an equivalence of data
$(M,K,\nnn^K,C,S)$ and $(H\to\P^1\times M,\nnn,P)$.
For slightly richer structures, such equivalences are formulated
in \cite[ch. VII]{S}, \cite[5.2]{H}, \cite[theorem 4.2]{HM}.

(ii) The main unfolding result theorem 2.5 in \cite{HM} applies also
to $(trTLP(0))$-structures, see \cite[remark 3.3 (vii)]{HM}.
In the situation in remark \ref{t4.3} (v), it applies, 
because its hypotheses are satisfied:
The generation condition (GC) above is a special case of the generation
condition in \cite[theorem 2.5]{HM}.
For any $x\in N$, the map 
\begin{eqnarray}\label{4.26}
T_xN\to K_x,\quad Y\mapsto C_Y\zeta,
\end{eqnarray}
is injective, so it satisfies the injectivity condition in 
\cite[theorem 2.5]{HM}. The unfolding result reads in our situation 
as follows.

For any $x\in N$, the germ of the tuple
$((N,x),K|_{(N,x)},\nnn^K,C,S)$ has a unique 
(up to isomorphism) unfolding 
to a tuple $((\www N,x),\www K,\nnn^{\www K},\www C,\www S)$
with the properties: $(\www N,x)\supset(N,x)$ is the germ at $x$ of
a manifold of dimension $\mu:=\rank K$, $\www K\to (\www N,x)$
is a vector bundle of rank $\mu$, $\nnn^{\www K}$ is a flat
connection on it, $\www C$ is a Higgs field on it
with $\nnn^{\www K}(\www C)=0$, and $\www S$ is a 
symmetric nondegenerate $\nnn^{\www K}$-flat and
Higgs field invariant pairing on $\www K$,
and finally, the map 
\begin{eqnarray}\label{4.27}
T_x\www N\to \www K_x,\quad Y\mapsto \www C_Y\zeta,
\end{eqnarray}
is an isomorphism.
On $(N,x)\subset (\www N,x)$ the tuple restricts to
the tuple $((N,x),K|_{(N,x)},\nnn^K,C,S)$.

The Higgs field endomorphisms $\www C_X,X\in\TT_{\www N}$,
induce a bundle $\www A$ of $\mu$-dimensional
commutative subalgebras $A_z\subset\textup{End}(K_z)$
for $z\in\www N$, and the map
\begin{eqnarray}\label{4.28}
T_z\www N\to A_z,\quad Y\mapsto C_Y
\end{eqnarray}
is an isomorphism. It induces a multiplication on
$\TT_{\www N}$ which turns out to give $\TT_{\www N}$
the structure of an $F$-manifold
(\cite[lemma 4.1 and lemma 4.3]{H}).

\medskip
(iii) We continue with the situation in (ii).
Choose an extension of the section $\zeta$ to a section
$\www\zeta$ in the bundle $\www K\to (\www N,x)$.
The isomorphism \eqref{4.25} extends to an isomorphism
\begin{eqnarray}\label{4.29}
\TT_{\www N}\to \OO(\www K),\quad Y\mapsto \www C_Y\www\zeta.
\end{eqnarray}
The section $\www\zeta$, the flat connection $\nnn^{\www K}$, 
the Higgs field $\www C$ and the pairing $\www S$
can be shifted to $\TT_{\www N}$ with this isomorphism.
The induced Higgs field gives the multiplication above:
That does not depend on the choice of $\www\zeta$.
But the other induced data depend on it.

One wishes an extension $\www\zeta$ such that the induced
connection on $\TT_{\www N}$ is torsion free.
One wishes a natural way to extend the notion of
a Frobenius like structure to the bundle $\www K\to\www N$.
And one wishes an extension $\www\zeta$ which
shifts this extension in the best possible way to $\TT_{\www N}$.

\medskip
(iv) The following special case gives rise to
Frobenius manifolds without Euler fields.
Consider a Frobenius like structure
$(M,K,\nnn^K,C,S,\zeta,V,(v_1,...,v_n))$ of order
$(n,1,2)$ with nondegenerate pairing $S$, 
$\nnn^K$-flat section $\zeta$, the conditions
(GC) and (wIC) and all $v_j\neq 0$ in the 1-dimensional
space $V$. Then the sections $C_{\paa_i}\zeta$ generate
the bundle $K$ and are $\nnn^K$-flat, the map 
$T_xM\to K_x,\quad Y\mapsto C_Y\zeta,$ is surjective with
1-dimensional kernel, $\mu=n-1$, 
the map \eqref{4.24} is an isomorphism,
the tuple $((N,x),K|_{(N,x)},\nnn^K,C,S)$ is its own
universal unfolding, and $N=\www N$.
Here $N$ becomes a Frobenius manifold (without Euler field).
The induced connection on $\TT_N$ is the one of the 
affine linear structure on $N$ and is torsion free.
It is also the Levi-Civita connection of the
metric on $\TT_M$ which is induced by $S$.
The restriction of the potential $L$ to $N$
is the potential of the Frobenius manifold.
\end{remarks}

\end{document}